\documentclass[a4paper,oneside,10pt]{article}%
\usepackage{amsmath}
\usepackage{amsfonts}
\usepackage{amssymb}
\usepackage{graphicx}
\usepackage{color}
\usepackage[square,numbers,sort&compress]{natbib}%
\providecommand{\U}[1]{\protect\rule{.1in}{.1in}}
\newtheorem{theorem}{Theorem}[section]

\newtheorem{corollary}[theorem]{Corollary}

\newtheorem{lemma}[theorem]{Lemma}

\numberwithin{equation}{section}

\begin{document}

\title{A filtering problem with uncertainty in observation}
\author{Shaolin Ji\thanks{Zhongtai Institute of Finance, Shandong University, Jinan,
Shandong 250100, PR China. jsl@sdu.edu.cn. This research is supported by
National Natural Science Foundation of China (No. 11571203), the Programme of
Introducing Talents of Discipline to Universities of China (No. B12023).}
\quad Chuiliu Kong\thanks{Corresponding author. Zhongtai Institute of Finance,
Shandong University, Jinan, Shandong 250100, PR China.
kclsdmath@mail.sdu.edu.cn.}\quad Chuanfeng Sun\thanks{School of Mathematical
Sciences, University of Jinan, Jinan, Shandong 250022, P.R. China.
sms\_suncf@ujn.edu.cn. This research is partially supported by the National
Natural Science Foundation of China (No. 11701214), the Natural Science
Foundation of Shandong Province (No. ZR2017BA032).}}
\date{}
\maketitle

\textbf{Abstract}. This paper is concerned with a generalized Kalman-Bucy
filtering model and corresponding robust problem under model uncertainty. We
find that this robust problem is equivalent to considering an estimate problem
under some sublinear operator. Therefore, we turn to obtaining the minimum
mean square estimator under a sublinear operator. By Girsanov theorem and
minimax theorem, we obtain the optimal estimator $\hat{x}_{t}$ of the signal
process $x_{t}$ for given time $t\in\lbrack0,T]$.


\section{Introduction}

Let $(\Omega,\mathcal{F},\{\mathcal{F}_{t}\}_{0\leq t\leq T},\emph{P})$ be a
complete filtered probability space equipped with a natural filtration
$\mathcal{F}_{t}=\sigma\{w(s),v(s);0\leq s\leq t\}$, $\mathcal{F}%
={\mathcal{F}_{T}}$, where $(w(\cdot),v(\cdot))$ is 2-dimensional standard
Brownian motion defined on the space, $T>0$ is a fixed real number. Suppose
that the signal process $({x}_{t})$ and the observation process $({m}_{t})$
under probability measure $P$ satisfy respectively
\begin{equation}
\left\{
\begin{array}
[c]{rl}%
dx_{t} & =(F_{t}x_{t}+f_{t})dt+dw_{t},\\
x(0) & =x_{0},\\
dm_{t} & =(G_{t}{x}_{t}+g_{t})dt+dv_{t},\\
m(0) & =0
\end{array}
\right.  \label{MO}%
\end{equation}
where the coefficients $F_{t},\ f_{t},\ G_{t},\ g_{t}$ are bounded, continuous
functions in $t$ and $x_{0}$ is a given constant. The classical Kalman-Bucy
filtering problem is to find the optimal estimator $\bar{x}_{t}$ such that
\begin{equation}
\min_{\zeta}E_{P}\Vert x_{t}-\zeta\Vert^{2}=E_{P}\Vert x_{t}-\bar{x}_{t}%
\Vert^{2}. \label{classical problem}%
\end{equation}

In 1961, Kalman and Bucy \cite{Kalman} gave the fundamental results of the
filtering problem which are the foundation of modern filtering theory (see
Bensoussan \cite{Bensoussan}, Liptser and Shiryaev \cite{Liptser} et al).
Based on the filtering technique, stochastic optimal control problems with
partial information (or observation) have been studied extensively. In the
field of finance and insurance, for example, Bensoussan and Keppo
\cite{Bensoussan-Keppo} and Lakner \cite{L} considered the optimal consumption
and portfolio investment problems of an investor who is interested in
maximizing his utilities from consumption and terminal wealth under partial
information; Xiong and Zhou \cite{Xiong-Zhou} considered the mean-variance
portfolio selection problems under partial information. In the field of
stochastic control, Duncan and Pasik-Dunan \cite{Duncan1} and \cite{Duncan2}
considered respectively the optimal control for a partially observed linear
stochastic system with an exponential quadratic cost and with fractional
brownian motions; Tang \cite{Tang} gave the maximum principle for partially
observed optimal control problems\ of stochastic differential equations; Wang
and Wu \cite{Wang-Wu} studied the Kalman-Bucy filtering equation of a certain
forward-backward stochastic differential equation system and solved a
partially observed linear quadratic optimal control problem, and so on. Some
fundamental researches based on forward-backward stochastic differential
equations are surveyed by Ma and Yong \cite{Ma-Yong} and Zhang \cite{Zhang}.

In 2002, Chen and Epstein \cite{Chen-Epstein} proposed a kind of model
uncertainty for continuous-time models which is the so called drift ambiguity.
Drift ambiguity models an agent's uncertainty about the drift of the
underlying Brownian motion. Moreover, in 2013, Epstein and Ji proposed more
general uncertainty models (see \cite{Epstein-Ji-1} and \cite{Epstein-Ji-2}
for details). In this paper, we introduce the following drift ambiguity in
\cite{Chen-Epstein} into model \eqref{MO} and focus on a corresponding robust
problem. Consider the generalized Kalman-Bucy filtering model under some
probability measure $P^{\theta}\in\mathcal{P}$:
\begin{equation}
\left\{
\begin{array}
[c]{rl}%
dx_{t} & =(F_{t}x_{t}+f_{t})dt+dw_{t},\\
x(0) & =x_{0},\\
dm_{t} & =(G_{t}{x}_{t}+g_{t}+\theta_{t})dt+dv_{t}^{\theta},\\
m(0) & =0
\end{array}
\right.  \label{generalize KB}%
\end{equation}
where $(w_{t})$ and $(v_{t}^{\theta})$ are Brownian motions under $P^{\theta}$
and the probability measure $P^{\theta}$ is regarded as an observer's
evaluation criterion for the signal process. Here the probability measure set
$\mathcal{P}$ denotes all the evaluation criterions by observers and
$\theta\in\Theta$ is called ambiguity parameter. Note that Ji, Li and Miao
\cite{Ji-Li-Miao} adopt a similar formulation in order to solve a dynamic
contract problem. Then, we naturally consider the following worst-case minimum
mean square estimate of the signal process $(x_{t})$:
\begin{equation}
\min_{\zeta}\sup_{P^{\theta}\in\mathcal{P}}E_{P^{\theta}}\Vert{x}_{t}%
-\zeta\Vert^{2}\label{intro-problem}%
\end{equation}
which is to minimize the maximum expected loss over a range of possible
models, an idea that goes back at least as far as Wald \cite{Wald} in 1945.
Allan and Cohen \cite{Allan-Cohen} studied this type of estimate problem under
nonlinear expectations by a control approach. Recently, Ji, Kong and Sun
\cite{Ji-Kong-Sun-2} considered a different generalized Kalman-Bucy filtering
model where the ambiguity parameters affect the evolution of signal process.

In fact, $\underset{P^{\theta}\in\mathcal{P}}{\sup}E_{P^{\theta}}[\cdot]$ can
be regarded as a sublinear operator $\mathcal{E}(\cdot)$ and the problem
\eqref{intro-problem} can be reformulated as a estimate problem under
sublinear operator:
\[
\min_{\zeta}\mathcal{E(}\Vert{x}_{t}-\zeta\Vert^{2}).
\]

The related literatures about the estimate problem under sublinear operators
include Sun and Ji \cite{Sun-Ji}, Ji, Kong and Sun \cite{Ji-Kong-Sun-1}. Sun
and Ji \cite{Sun-Ji} introduced a new conditional nonlinear expectation for
bounded random variables which is based on the minimum mean square estimator
for sublinear operators. However, the boundedness assumption for random
variables has great limitations. Therefore, Ji, Kong and Sun
\cite{Ji-Kong-Sun-1} deleted the boundedness assumption and generalized the
corresponding results to the case in which the random variables fall in the
space $L_{\mathcal{F}}^{2+\epsilon}(\Omega,P)$ where $\epsilon$ is a constant
such that $\epsilon\in(0,1)$.

Under some mild conditions, we prove that the optimal estimator $\hat{x}$\ and
the optimal probability measure $P^{\theta^{\ast}}$ exist. It results that we
only need to consider the classical Kalman-Bucy filtering problem under the
probability measure $P^{\theta^{\ast}}$. Moreover, in some special cases, the
optimal estimator $\hat{x}$ can be decomposed to two parts. One part is the
optimal estimator of the signal process under the probability measure $P$ and
the other part contains the parameter $\theta^{\ast}$ (see Corollary
\ref{decomposition} for details).

The remainder of the paper proceeds as follows. In section 2, after a brief
recall of the Kalman-Bucy filter and the drift ambiguity, a generalized robust
Kalman-Bucy filtering problem is introduced. In section 3, the main general
results are given and we consider a special case to further explain our results.

\section{Problem formulation}

Let $w(\cdot)$ and $v(\cdot)$ be $n$-dimensional and $m$-dimensional
independent Brownian motions defined on a complete filtered probability space
$(\Omega,\mathcal{F},\{\mathcal{F}_{t}\}_{0\leq t\leq T},\emph{P})$ where
$\mathcal{F}_{t}=\sigma\{w(s),v(s);0\leq s\leq t\}$, $\mathcal{F}%
={\mathcal{F}_{T}}$ and $T>0$ be a fixed terminal time. The means of
$w(\cdot)$ and $v(\cdot)$ are zero and the covariance matrices are $Q(\cdot)$
and $R(\cdot)$ respectively. The matrix $R(\cdot)$ is uniformly positive
definite. Denote by $\mathbb{R}^{n}$ the $n$-dimensional real Euclidean space
and $\mathbb{R}^{n\times k}$ the set of $n\times k$ real matrices. Let
$\langle\cdot,\cdot\rangle$ (resp. $\Vert\cdot\Vert$) denote the usual scalar
product (resp. usual norm) of $\mathbb{R}^{n}$ and $\mathbb{R}^{n\times k}$.
The scalar product (resp. norm) of $M=(m_{ij})$, $N=(n_{ij})\in\mathbb{R}%
^{n\times k}$ is denoted by $\langle M,N\rangle=tr\{MN^{\intercal}\}$ (resp.
$\Vert M\Vert=\sqrt{\langle M,M\rangle}$), where the superscript $^{\intercal
}$ denotes the transpose of vectors or matrices. For a $\mathbb{R}^{n}$-valued
vector $x=(x_{1},\cdot\cdot\cdot,x_{n})^{\intercal}$, $|x|:=(|x_{1}%
|,\cdot\cdot\cdot,|x_{n}|)^{\intercal}$; for two $\mathbb{R}^{n}$-valued
vectors $x$ and $y$, $x\leq y$ means that $x_{i}\leq y_{i}$ for $i=1,\cdot
\cdot\cdot,n$.

Through out this paper, $0$ denotes the matrix/vector with appropriate
dimension whose all entries are zero and $\epsilon$ is a constant such that
$0<\epsilon<1$.

Suppose that the signal process $({x}_{t})\in\mathbb{R}^{n}$ and the
observation process $({m}_{t})\in\mathbb{R}^{m}$ under probability measure $P$
satisfy model \eqref{MO} where $F_{t}\in\mathbb{R}^{n\times n},\ G_{t}%
\in\mathbb{R}^{m\times n},\ f_{t}\in\mathbb{R}^{n},\ g_{t}\in\mathbb{R}^{m}$
are bounded, continuous functions in $t$, $x_{0}\in\mathbb{R}^{n}$ is a given
constant vector. Let the filtration $\mathcal{Z}_{t}=\sigma\{m(s);0\leq s\leq
t\}$ be the set of observable events up to time $t$. By the Kalman-Bucy
filtering theory (see Bensoussan \cite{Bensoussan}, Kalman and Bucy
\cite{Kalman} and Liptser and Shiryaev \cite{Liptser} et al), the optimal
solution $\bar{x}_{t}=E_{P}(x_{t}|\mathcal{Z}_{t})$ of problem
\eqref{classical problem} is governed by
\begin{equation}
\left\{
\begin{array}
[c]{ll}%
d\bar{x}_{t} & =(F_{t}\bar{x}_{t}+f_{t})dt+P_{t}G_{t}^{\intercal}R_{t}%
^{-1}dI_{t},\\
\bar{x}(0) & =x_{0},
\end{array}
\right.  \label{classical Kalman}%
\end{equation}
and the variance of estimate error $P_{t}=E_{P}[(x_{t}-\bar{x}_{t})(x_{t}%
-\bar{x}_{t})^{\intercal}]$ is governed by
\begin{equation}
\left\{
\begin{array}
[c]{rl}
& \frac{dP_{t}}{dt}=F_{t}P_{t}+P_{t}F_{t}^{\intercal}-P_{t}G_{t}^{\intercal
}R_{t}^{-1}G_{t}P_{t}+Q_{t},\\
& P(0)=0
\end{array}
\right.  \label{Riccati}%
\end{equation}
where $I_{t}=m_{t}-\int_{0}^{t}(G_{s}\bar{x}_{s}+g_{s})ds$ is called
innovation process under probability measure $P$ which is a Wiener process
adapted to $\{\mathcal{Z}_{t}\}$. Furthermore, the filtration $\mathcal{I}%
_{t}=\sigma\{I(s);0\leq s\leq t\}$ equals to $\mathcal{Z}_{t}$ for any time
$t\in\lbrack0,T]$.

Now we are ready to give the drift ambiguity model. For a fixed $\mathbb{R}%
^{m}$-valued nonnegative constant vector $\mu$, denote by $\Theta$ the set of
all the $\mathbb{R}^{m}$-valued progressively measurable processes
$(\theta_{t})$ with $|\theta_{t}|\leq\mu$. Define
\begin{equation}
\mathcal{P}=\{P^{\theta}\big|\frac{dP^{\theta}}{dP}=f_{T}^{P^{\theta}%
}\ \mbox{with}\ \theta\in\Theta\} \label{CD}%
\end{equation}
where
\[
f_{T}^{P^{\theta}}=\frac{dP^{\theta}}{dP}=\exp(\int_{0}^{T}\theta
_{t}^{\intercal}dv_{t}-\frac{1}{2}\int_{0}^{T}\Vert\theta_{t}\Vert^{2}dt).
\]

Due to the boundness of $\theta$, the Novikov's condition holds (see Karatzas
and Shreve \cite{K-S}). Therefore, $P^{\theta}$ defined by \eqref{CD} is a
probability measure which is equivalent to the probability measure $P$ and the
processes $(w_{t})$ and $(v_{t}^{\theta})$ where $v_{t}^{\theta}=v_{t}%
-\int_{0}^{t}\theta_{s}ds$ are Brownian motions under this probability measure
$P^{\theta}$ by Girsanov theorem. Then, with this generalized model
\eqref{generalize KB} under probability measure $P^{\theta}$, we consider
naturally the following robust problem:
\begin{equation}
\inf_{\zeta\in L_{\mathcal{Z}_{t}}^{2+\epsilon}(\Omega,P,\mathbb{R}^{n})}%
\sup_{P^{\theta}\in\mathcal{P}}E_{P^{\theta}}\Vert{x}_{t}-\zeta\Vert^{2},
\label{robust problem}%
\end{equation}
where $L_{\mathcal{Z}_{t}}^{2+\epsilon}(\Omega,P,\mathbb{R}^{n})$ is the set
of all the $\mathbb{R}^{n}$-valued $(2+\epsilon)$ integrable $\mathcal{Z}_{t}%
$-measurable random variables.

However, if we denote $\mathcal{E}(\cdot)=\underset{P^{\theta}\in
\mathcal{P}}{\sup}E_{P^{\theta}}[\cdot]$ which can be regarded as a sublinear
operator, then the above robust problem can be considered as an estimate
problem of the signal process under this sublinear operator $\mathcal{E}%
(\cdot)$. In more details, given the observation information $\{\mathcal{Z}%
_{t}\}$, we intend to find the optimal estimator $\hat{x}_{t}$ of the signal
process $(x_{t})$ at time $t\in\lbrack0,T]$ such that
\begin{equation}
\mathcal{E}\Vert{x}_{t}-\hat{x}_{t}\Vert^{2}=\inf_{\zeta\in\mathcal{K}_{t}%
}\mathcal{E}\Vert{x}_{t}-\zeta\Vert^{2}, \label{OP}%
\end{equation}
where
\[
\mathcal{K}_{t}=\{\zeta:\Omega\rightarrow\mathbb{R}^{n};\ \zeta\in
L_{\mathcal{Z}_{t}}^{2+\epsilon}(\Omega,P,\mathbb{R}^{n})\}.
\]
\textbf{Remark 2.1.} The optimal solution $\hat{x}_{t}$ of problem \eqref{OP}
is called minimum mean square estimator. It is also regarded as a minimax
estimator in statistical decision theory. If the sublinear operator
$\mathcal{E}(\cdot)$ degenerates to linear expectation operator, then
$\mathcal{P}^{\theta}$ contains only one probability measure $P$. In this
case, it is well known that the minimum mean square estimator $\hat{x}_{t}$ is
just the conditional expectation $E_{P}({x}_{t}|\mathcal{Z}_{t})$.

\section{Main results}

In this section, we study the minimum mean square estimator $\hat{{x}}_{t}$ of
problem \eqref{OP} for some time $t\in[0,T]$. Without loss of generality, we
only prove one dimensional case and the multidimensional case can be proved similarly.

\begin{lemma}
\label{LW} The set $\{\frac{dP^{\theta}}{dP}:P^{\theta}\in\mathcal{P}\}\subset
L^{1+\frac{2}{\epsilon}}(\Omega,\mathcal{F},P)$ is $\sigma(L^{1+\frac
{2}{\epsilon}}(\Omega,\mathcal{F},P),L^{1+\frac{\epsilon}{2}}(\Omega
,\mathcal{F},P))$-compact and the set $\mathcal{P}$ is convex.
\end{lemma}

$\mathbf{Proof}.$ By Lemma 1 in Girsanov \cite{Girsanov} and the boundness of
$\theta$, the set $\{\frac{dP^{\theta}}{dP}:P^{\theta}\in\mathcal{P}\}\subset
L^{1+\frac{2}{\epsilon}}(\Omega,\mathcal{F},P)$ space. According to Simons
\cite{Simons}, Chapter 1, Theorem 4.1, the set $\{\frac{dP^{\theta}}%
{dP}:P^{\theta}\in\mathcal{P}\}$ is $\sigma(L^{1+\frac{2}{\epsilon}}%
(\Omega,\mathcal{F},P), L^{1+\frac{\epsilon}{2}}\newline(\Omega,\mathcal{F}%
,P))$-compact.

The set $\mathcal{P}$ is convex which\ can be referred to Chen and Epstein
\cite{Chen-Epstein}. Let $\theta_{1}$ and $\theta_{2}$ belong to the set
$\Theta$. $f^{P^{\theta_{i}}},\ i=1,2$ denote the exponential martingales
respectively with
\[
f_{t}^{P^{\theta_{i}}}=\exp(\int_{0}^{t}\theta_{i,s}dv_{s}-\frac{1}{2}\int%
_{0}^{t}\theta_{i,s}^{2}ds)
\]
and
\[
df_{t}^{P^{\theta_{i}}}=f_{t}^{P^{\theta_{i}}}\theta_{i,t}dv_{t}.
\]
Let $0\leq\lambda_{i}\leq1,\ i=1,2$ be constants with $\lambda_{1}+\lambda
_{2}=1$ and
\[%
\begin{array}
[c]{rl}%
\theta_{t}^{\lambda} & =\frac{\lambda_{1}\theta_{1,t}f_{t}^{P^{\theta_{1}}%
}+\lambda_{2}\theta_{2,t}f_{t}^{P^{\theta_{2}}}}{\lambda_{1}f_{t}%
^{P^{\theta_{1}}}+\lambda_{2}f_{t}^{P^{\theta_{2}}}}.
\end{array}
\]
Since $f^{P^{\theta_{i}}}>0,\ i=1,2$, the process $(\theta_{t}^{\lambda})$
belongs to the set $\Theta$, which implies that the set $\Theta$ is
stochastically convex. Moreover, it is also easy to calculate that
\[
d(\lambda_{1}f_{t}^{P^{\theta_{1}}}+\lambda_{2}f_{t}^{P^{\theta_{2}}%
})=(\lambda_{1}f_{t}^{P^{\theta_{1}}}+\lambda_{2}f_{t}^{P^{\theta_{2}}}%
)\theta_{t}^{\lambda}dv_{t}.
\]
Therefore, the set $\mathcal{P}$ is convex.\newline\rightline{$\square$}
\textbf{Remark 3.1.} By Lemma \ref{LW}, Lemma 1 in \cite{Girsanov} and Theorem
6.3 in Chapter 1 of \cite{Yong-Zhou}, the signal process $(x_{t})$ is
$(4+2\epsilon)$ integrable, the set $\{\frac{dP^{\theta}}{dP}:P^{\theta}%
\in\mathcal{P}\}$ is uniformly normed bounded in $L^{1+\frac{2}{\epsilon}%
}(\Omega,\mathcal{F},P)$ space and also $\sigma(L^{1+\frac{2}{\epsilon}%
}(\Omega,\mathcal{F},P),L^{1+\frac{\epsilon}{2}}(\Omega,\mathcal{F}%
,P))$-compact. Therefore, we can apply the results in Ji, Kong and Sun
\cite{Ji-Kong-Sun-1} to guarantee that the optimal solution of problem
\eqref{OP} exists.

By Lemma \ref{LW}, we can apply the minimax theorem (see Theorem B.1.2 in Pham
\cite{Pham}) to problem \eqref{robust problem} which leads to the following theorem.

\begin{theorem}
\label{ER} For a given $t\in\lbrack0,T]$, there exists a $\theta^{\ast}%
\in\Theta$ such that
\[
\label{ERE}\inf_{\zeta\in\mathcal{K}_{t}}\mathcal{E}\Vert{x}_{t}-\zeta
\Vert^{2}=\inf_{\zeta\in\mathcal{K}_{t}}\sup_{P^{\theta}\in\mathcal{P}%
}E_{P^{\theta}}\Vert{x}_{t}-\zeta\Vert^{2}=\inf_{\zeta\in\mathcal{K}_{t}%
}E_{P^{\theta^{\ast}}}\Vert{x}_{t}-\zeta\Vert^{2}.
\]

\end{theorem}

\textbf{Proof}. Denote $f_{n}=\frac{dP^{\theta_{n}}}{dP}$ and choose a
sequence $\{f_{n}\}_{n\geq1}$ such that
\[
\lim_{n\rightarrow\infty}\inf_{\zeta\in\mathcal{K}_{t}}E_{P}[f_{n}({x}%
_{t}-\zeta)^{2}]=\lim_{n\rightarrow\infty}\inf_{\zeta\in\mathcal{K}_{t}%
}E_{P^{\theta_{n}}}[({x}_{t}-\zeta)^{2}]=\sup_{P^{\theta}\in\mathcal{P}}%
\inf_{\zeta\in\mathcal{K}_{t}} E_{P^{\theta}}[({x}_{t}-\zeta)^{2}].
\]
By Koml\'{o}s theorem in Pham \cite{Pham}, we know that there exist a
subsequence $\{f_{n_{k}}\}_{k\geq1}\subset\{f_{n} \}_{n\geq1}$ and $f^{\ast
}\in L^{1}(\Omega,\mathcal{F},P)$ space such that
\[
\lim_{m\rightarrow\infty}\dfrac{1}{m}\sum_{k=1}^{m}f_{n_{k}}=f^{\ast
},\ P-a.s..
\]
Let $g_{m}=\dfrac{1}{m}\sum_{k=1}^{m}f_{n_{k}}$. We have $g_{m}%
\xrightarrow{P-a.s.}f^{\ast}$ and
\begin{equation}%
\begin{array}
[c]{rl}
& \sup_{P^{\theta}\in\mathcal{P}}\inf_{\zeta\in\mathcal{K}_{t}} E_{P^{\theta}%
}[({x}_{t}-\zeta)^{2}]=\lim_{n\rightarrow\infty}\inf_{\zeta\in\mathcal{K}_{t}%
}E_{P^{\theta_{n}}}[({x}_{t}-\zeta)^{2}]=\lim_{k\rightarrow\infty}\inf
_{\zeta\in\mathcal{K}_{t}}E_{P^{\theta_{n_{k}}}}[({x}_{t}-\zeta)^{2}]\\
& =\lim_{m\rightarrow\infty}\frac{1}{m}\sum_{k=1}^{m}\inf_{\zeta\in
\mathcal{K}_{t}}E_{P^{\theta_{n_{k}}}}[({x}_{t}-\zeta)^{2}]\leq\lim
_{m\rightarrow\infty}\inf_{\zeta\in\mathcal{K}_{t}}\frac{1}{m}\sum_{k=1}%
^{m}E_{P^{\theta_{n_{k}}}}[({x}_{t}-\zeta)^{2}]\\
& =\lim_{m\rightarrow\infty}\inf_{\zeta\in\mathcal{K}_{t}}E_{P}[g_{m}({x}%
_{t}-\zeta)^{2}].
\end{array}
\end{equation}

By Lemma 1 in \cite{Girsanov}, for any given constants $p>1$ and $m$, we have
$\ E_{P}(g_{m})^{K}\leq M$ \ where $K=(1+\frac{2}{\epsilon})p$ and
$M=\exp(\frac{K^{2}-K}{2}\mu^{2}T)$. Then, we have $\left\{  |g_{m}%
|^{1+\frac{2}{\varepsilon}}:m=1,2,\cdot\cdot\cdot\right\}  $ is uniformly
integrable. Therefore, $\ g_{m}%
\xrightarrow{L^{1+\frac{2}{\epsilon}}(\Omega, \mathcal{F}, P)}f^{\ast}$ and
$f^{\ast}\in L^{1+\frac{2}{\epsilon}}(\Omega,\mathcal{F},P)$. According to the
convexity weak compactness of set $\{\frac{dP^{\theta}}{dP}: P^{\theta}%
\in\mathcal{P}\}$, there exists a $\theta^{*}$ such that $\frac{dP^{\theta
^{*}}}{dP}=f^{*}$ and the following relations hold
\[%
\begin{array}
[c]{rl}%
\sup_{P^{\theta}\in\mathcal{P}}\inf_{\zeta\in\mathcal{K}_{t}}E_{P^{\theta}%
}[({x}_{t}-\zeta)^{2}] & \geq\inf_{\zeta\in\mathcal{K}_{t}}E_{P^{\theta^{*}}%
}[({x}_{t}-\zeta)^{2}]\\
& =\inf_{\zeta\in\mathcal{K}_{t}}E_{P}[f^{*}({x}_{t}-\zeta)^{2}]\\
& =\inf_{\zeta\in\mathcal{K}_{t}}E_{P}[\lim_{m\rightarrow\infty} g_{m}({x}%
_{t}-\zeta)^{2}]\\
& \geq\limsup_{m\rightarrow\infty}\inf_{\zeta\in\mathcal{K}_{t}}E_{P}%
[g_{m}({x}_{t}-\zeta)^{2}]\\
& \geq\sup_{P^{\theta}\in\mathcal{P}}\inf_{\zeta\in\mathcal{K}_{t}%
}E_{P^{\theta}}[({x}_{t}-\zeta)^{2}]
\end{array}
\]
where the second $^{\prime}\geq^{\prime}$ is based on upper semi-continuous
property. It follows that
\[
\sup_{P^{\theta}\in\mathcal{P}}\inf_{\zeta\in\mathcal{K}_{t}}E_{P^{\theta}%
}[({x}_{t}-\zeta)^{2}]=\inf_{\zeta\in\mathcal{K}_{t}}E_{P^{\theta^{*}}}%
[({x}_{t}-\zeta)^{2}].
\]
By minimax theorem, we obtain
\[
\sup_{P^{\theta}\in\mathcal{P}}\inf_{\zeta\in\mathcal{K}_{t}}E_{P^{\theta}%
}[({x}_{t}-\zeta)^{2}]=\inf_{\zeta\in\mathcal{K}_{t}}\sup_{P^{\theta}%
\in\mathcal{P}}E_{P^{\theta}}[({x}_{t}-\zeta)^{2}].
\]
It results that
\[
\inf_{\zeta\in\mathcal{K}_{t}}\sup_{P^{\theta}\in\mathcal{P}}E_{P^{\theta}%
}[({x}_{t}-\zeta)^{2}]=\inf_{\zeta\in\mathcal{K}_{t}}E_{P^{\theta^{*}}}%
[({x}_{t}-\zeta)^{2}].
\]
\rightline{$\square$}

Once we find the optimal $\theta^{\ast}$, model \eqref{generalize KB} and
problem \eqref{robust problem} can be expressed under the new probability
measure $P^{\theta^{\ast}}$ correspondingly. In more details, for this
filtered probability space $(\Omega,\mathcal{F},\{\mathcal{F}_{t}\}_{0\leq
t\leq T},P^{\theta^{\ast}})$, the processes $(x_{t})$ and $(m_{t})$ satisfy
respectively
\begin{equation}
\left\{
\begin{array}
[c]{rl}%
{dx}_{t} & =(F_{t}{x}_{t}+f_{t})dt+dw_{t},\\
{x}(0) & =x_{0},\\
{dm}_{t} & =(G_{t}x_{t}+g_{t}+\theta_{t}^{\ast})dt+dv_{t}^{\theta^{\ast}},\\
{m}(0) & =0
\end{array}
\right.  \label{MO1}%
\end{equation}
and problem \eqref{robust problem} turns into the minimum mean square estimate
problem under probability measure $P^{\theta^{\ast}}$:%
\begin{equation}
E_{P^{\theta^{\ast}}}\Vert{x}_{t}-\hat{x}_{t}\Vert^{2}=\inf_{\zeta
\in\mathcal{K}_{t}}E_{P^{\theta^{\ast}}}\Vert{x}_{t}-\zeta\Vert^{2}.
\label{NP2}%
\end{equation}

With the above theorem, we consider the following estimate problem which is a
Kalman-Bucy filtering problem with the parameter $\theta^{\ast}$:
\begin{equation}
E_{P^{\theta^{\ast}}}\Vert{x}_{t}-\hat{\zeta}\Vert^{2}=\inf_{\zeta\in
\bar{\mathcal{K}}_{t}}E_{P^{\theta^{\ast}}}\Vert{x}_{t}-\zeta\Vert^{2}
\label{KB}%
\end{equation}
where $\bar{\mathcal{K}}_{t}=\{\zeta:\Omega\rightarrow\mathbb{R}^{n}%
;\ \zeta\in L_{\mathcal{Z}_{t}}^{2}(\Omega,P^{\theta^{\ast}},\mathbb{R}%
^{n})\}$.

The model \eqref{MO1} and problem \eqref{KB} constitute a classical
construction for a linear, partially observable system with a parameter
$\theta^{\ast}$. This estimate problem is to characterize the conditional
distribution $P^{\theta^{\ast}}({x}_{t}\in A|\mathcal{Z}_{t})$, where $A$ is a
Borel set in $\mathbb{R}^{n}$. Then, we are in the realm of Kalman-Bucy
filtering and it is well known(see \cite{Kalman} and \cite{Liptser}) that the
conditional distribution is again Gaussian and conditional mean $\hat{x}%
_{t}=E_{P^{\theta^{\ast}}}({x}_{t}|\mathcal{Z}_{t})$ solves the following
equation:%
\begin{equation}
\left\{
\begin{array}
[c]{rl}%
d\hat{x}_{t} & =(F_{t}\hat{x}_{t}+f_{t})dt+(P_{t}G^{\intercal}_{t}%
+\widehat{x_{t}\theta^{\ast\intercal}_{t}}-\hat{x}_{t}\widehat{\theta
^{\ast\intercal}_{t}})R_{t}^{-1}d\hat{I}_{t},\\
\hat{x}(0) & =x_{0}%
\end{array}
\right.  \label{KBequation}%
\end{equation}
where $\widehat{\theta_{t}^{\ast}}=E_{P^{\theta^{\ast}}}[\theta^{\ast}%
_{t}|\mathcal{Z}_{t}]$, $\widehat{x_{t}\theta^{\ast\intercal}_{t}%
}=E_{P^{\theta^{\ast}}}[x_{t}\theta^{\ast\intercal}_{t}|\mathcal{Z}_{t}]$,
$\hat{I}_{t}={m}_{t}-\int_{0}^{t}(G_{s}\hat{x}_{s}+g_{s}+\widehat{\theta
^{*}_{s}})ds$ is $\mathcal{Z}_{t}$-measurable Brownian motion and the variance
of error equation $P_{t}=E_{P^{\theta^{\ast}}}[(x_{t}-\hat{x}_{t}%
)^{2}|\mathcal{Z}_{t}]=E_{P^{\theta^{\ast}}}[(x_{t}-\hat{x}_{t})^{2}]$
satisfies
\begin{equation}
\left\{
\begin{array}
[c]{rl}%
\frac{dP_{t}}{dt} & =F_{t}P_{t}+P_{t}F^{\intercal}_{t}-E_{P^{\theta^{\ast}}%
}[(P_{t}G^{\intercal}_{t}+\widehat{x_{t}\theta^{\ast\intercal}_{t}}-\hat
{x}_{t}\widehat{\theta^{\ast\intercal}_{t}})R^{-1}_{t}(G_{t}P_{t}%
+\widehat{\theta^{\ast}_{t}x^{\intercal}_{t}}-\widehat{\theta^{\ast}_{t}}%
\hat{x^{\intercal}}_{t})]+Q_{t},\\
P(0) & =0.
\end{array}
\right.  \label{NewRiccati}%
\end{equation}

So far, the optimal estimator of problem \eqref{KB} has been obtained. Next,
we expound that this solution $\hat{x}_{t}$ is also the optimal estimator of
problem \eqref{OP} at time $t\in[0,T]$.

\begin{theorem}
\label{ER1} Under the above assumptions, $\hat{x}_{t}$ governed by equation
\eqref{KBequation} is also the optimal solution of problem \eqref{OP} for any
time $t\in[0,T]$.
\end{theorem}

\textbf{Proof.} Note that
\begin{equation}
\inf_{\zeta\in\mathcal{K}_{t}}\sup_{P^{\theta}\in\mathcal{P}}E_{P^{\theta}%
}({x}_{t}-\zeta)^{2}=\inf_{\zeta\in\mathcal{K}_{t}}E_{P^{\theta^{\ast}}}%
({x}_{t}-\zeta)^{2}\geq\inf_{\zeta\in\bar{\mathcal{K}}_{t}}E_{P^{\theta^{\ast
}}}({x}_{t}-\zeta)^{2}.
\end{equation}

In addition, since $F_{t},\ G_{t},\ f_{t}\ \mbox{and}\ g_{t}$ are bounded
continuous functions in $t$ and $\theta^{\ast}$ is bounded, it is easy to
verify that $\hat{x}_{t}$ is not only square integrable but also
$(4+2\epsilon)$ integrable under probability measure $P^{\theta^{\ast}}$ by
Theorem 6.3 in Chapter 1 of \cite{Yong-Zhou}. Then, the solution $\hat{x}_{t}$
of equation \eqref{KBequation} also belongs to $\mathcal{K}_{t}$. It yields
that $\hat{x}_{t}$\ is the optimal solution of problem \eqref{OP} at time
$t\in\lbrack0.T]$.\newline\rightline{$\square$}

\begin{corollary}
If the optimal $\theta^{\ast}_{t}$ adapted to subfiltration $\mathcal{Z}_{t}$,
then the optimal estimator $\hat{x}_{t}$ satisfies the following simpler
equation.
\begin{equation}
\left\{
\begin{array}
[c]{rl}%
d\hat{x}_{t} & =(F_{t}\hat{x}_{t}+f_{t})dt+P_{t}G^{\intercal}_{t}R_{t}%
^{-1}d\hat{I}_{1,t},\\
\hat{x}(0) & =x_{0}%
\end{array}
\right.  \label{SKBequation}%
\end{equation}
where $\hat{I}_{1,t}={m}_{t}-\int_{0}^{t}(G_{s}\hat{x}_{s}+g_{s}+\theta
^{*}_{s})ds$ is $\mathcal{Z}_{t}$-measurable Brownian motion and $P_{t}$
reduces to equation \eqref{Riccati}.
\end{corollary}

Define $A(t,s)=P_{s}G_{s}R_{s}^{-1}\exp^{\int_{s}^{t}(F_{r}-P_{r}G_{r}%
^{2}R_{r}^{-1})dr}$, which is the impulse response of the classical
Kalman-Bucy filter. After some simple calculations, the optimal estimator
$\hat{x}$ can be decomposed to two parts. One part is the optimal estimator of
the signal process under the probability measure $P$ and the other part
contains the parameter $\theta^{\ast}$.

\begin{corollary}
\label{decomposition} If the optimal $\theta^{\ast}_{t}$ adapted to
subfiltration $\mathcal{Z}_{t}$, with equations \eqref{classical Kalman} and
\eqref{SKBequation}, then the optimal estimator $\hat{x}_{t}$ for any time
$t\in\lbrack0.T]$ can be expressed as
\begin{equation}
\hat{x}_{t}=\bar{x}_{t}-\int_{0}^{t}A(t,s)\theta_{s}^{\ast}ds.
\label{expression2}%
\end{equation}
where $\bar{x}_{t}$ is defined by equation \eqref{classical Kalman}.
\end{corollary}

Similar to Theorem 5.6 in Sun and Ji \cite{Sun-Ji}, we give a sufficient and
necessary condition for the existence of the optimal estimator in the
following corollary.

\begin{corollary}
For a given $t\in[0,T]$, $\hat{x}_{t}$ is the optimal solution of problem
\eqref{OP} if and only if it is the solution of the following equation
\begin{equation}
\inf_{\zeta\in\mathcal{K}_{t}}\mathcal{E}[(x_{t}-\hat{x}_{t})(x_{t}%
-\zeta)]=\mathcal{E}(x_{t}-\hat{x}_{t})^{2}.
\end{equation}

\end{corollary}

\end{document}